\magnification=1200
\documentstyle{amsppt}
\document
\centerline{$A_\infty$-STRUCTURES AND DIFFERENTIALS OF THE}
\centerline{ADAMS SPECTRAL SEQUENCE}
\vskip 6pt
\centerline{\sl V.A.Smirnov}
\vskip .5 cm
The Adams spectral sequence was invented by J.F.Adams [1] almost fifty
years  ago  for  calculations of stable homotopy groups of topological
spaces and in particular of spheres.  The calculation of differentials
of  this  spectral  sequence  is  one of the most difficult problem of
Algebraic Topology. Here we consider an approach to solve this problem
in   the  case  of $\Bbb Z/2$ coefficients and find inductive formulas
for the differentials. It is  based  on  the $A_\infty$-structures [2],
operad methods [3], [4],  [5],  [6]  and functional homology operations
[7], [8], [9].

\vskip .5cm
\centerline{CONTENTS}
\vskip 6pt
1. The Bousfield-Kan spectral sequences

2. $E_\infty$-algebras and $E_\infty$-coalgebras

3. $A_\infty$-cosimplicial objects

4. $E_\infty$-structure on the Bousfield-Kan spectral sequence

5. The  homology  of  an  $E_\infty$-operad and the Milnor coalgebra

6. Degenerating $A_\infty$-structures

7. Functorial homology operations

8. Homology operations for the operad $E$

9. $\cup_\infty-A_\infty$-Hopf algebras

10. The calculation of the differentials

\vskip .5cm
\centerline{1. The Bousfield-Kan spectral sequences}
\vskip 6pt
Consider the Bousfield-Kan spectral sequence [10], which is one of the
most general spectral sequence of the homotopy groups.

Let $R$ be a field. Given a simplicial set $Z$ denote by $RZ$ the free
simplicial $R$-module generated by  $Z$.  There  is  the  cosimplicial
resolution
$$R^*Z:\quad Z@>\delta^0>>RZ@>\delta^0,\delta^1>>R^2Z\to\dots\to
R^nZ@>\delta^0,\dots, \delta^n>>R^{n+1}Z\to\dots.$$
This resolution was used by Bousfield and Kan [10]  to  construct  the
spectral  sequence  of the homotopy groups of $Z$ with coefficients in
$R$.

The $E^1$-term of this spectral sequence is expressed by  the  complex
$$H_*(Z;R)\to  H_*(RZ;R)\to\dots\to  H_*(R^{n-1}Z;R)\to H_*(R^nZ;R)\to
\dots.$$ Higher differentials are expressed by the homology operations
$$d_m\colon  H_*(R^{n-1}Z;R)\to  H_*(R^{n+m-1}Z;R).$$ In [7],  [8] the
homology operations were defined as partial and multi-valued mappings.
However there is a general method to choose the homology operations to
be usual homomorphisms. The corresponding theory was developed in [9].
Recall main definitions.

For a chain complex $X$ denote by $X_*$ its homology, $X_*=H_*(X)$.
Fix chain mappings $\xi\colon  X_*\to  X$,  $\eta\colon  X\to  X_*$ and
chain homotopy $h\colon X\to X$ satisfying the relations
$$\eta\circ\xi=Id,\quad d(h)=\xi\circ\eta-Id,\quad   h\circ\xi=0,\quad
\eta\circ h=0,\quad h\circ h=0.$$
Consider a sequence of mappings
$$f^1\colon X^1\to X^2,~\dots,~f^n\colon X^n\to X^{n+1}$$  and  define
functional homology operations
$$H_*(f^n,\dots,f^1)\colon X^1_*\to        X^{n+1}_*$$
putting
$$H_*(f^n,\dots,f^1)=\eta\circ        f^n\circ        h\circ\dots\circ
f^1\circ\xi.$$
Direct calculations show that the following relations are satisfied
$$\gather\sum_{i=1}^{n-1}(-1)^{n-i+1}H_*(f^n,\dots,       f^{i+1}\circ
f^i,\dots,f^1)=\\
\sum_{i=1}^{n-1}(-1)^{n-i}H_*(f^n,\dots,f^{i+1})\circ
H_*(f^i,\dots,f^1).\endgather$$

Functional homology  operations  may  be  defined  not  only  for  the
category of chain complexes but in some other situations, for example,
for the category of simplicial modules.

Directly from the definition it follows that higher  differentials  of
the Bousfield-Kan  spectral  sequence  are expressed by the functional
homology operations
$$H_*(\delta,\dots,\delta)\colon                    H_*(R^{n-1}Z;R)\to
H_*(R^{n+m-1}Z;R).$$ So we have

{\bf Theorem 1.} {\sl The differentials of the  Bousfield-Kan  spectral
sequence are expressed by the functional homology operations
$$H_*(\delta,\dots,\delta)\colon               H_*(R^{n-1}Z,R)\to
H_*(R^{n+m-1}Z,R).$$ These operations determine on the $E^1$-term  new
differential.  The homology of the corresponding complex is isomorphic
to the $E^\infty$-term of the spectral sequence.}

As it was proved in [5] instead of the Bousfield-Kan cosimplicial object
we may consider the following cosimplicial object
$$F^*(C,RZ):\quad RZ@>\delta^0,\delta^1>> CRZ\to\dots\to C^{n-1}RZ
@>\delta^0, \dots,\delta^n>>C^nRZ\to\dots,  $$ where $CRZ$ is the free
commutative simplicial coalgebra generated by $RZ$.

The $E^1$-term  of the corresponding spectral sequence is expressed by
the complex
$$H_*(Z;R)=\pi_*(RZ)\to \pi_*(CRZ)\to\dots\to      \pi_*(C^{n-1}RZ)\to
\pi_*(C^nRZ)\to\dots.$$ Directly from the definition it  follows  that
higher differentials  of  this  spectral sequence are expressed by the
functional homology operations
$$H_*[\delta,\dots,\delta]\colon\pi_*(C^nRZ)\to\pi_*(C^{n+m}RZ).$$
Moreover there is a cosimplicial mapping
$$\CD RZ@>>>R^2Z@>>>\dots @>>>R^{n+1}Z@>>>\dots\\
@V=VV @VVV @. @VVV @.\\
RZ@>>>CRZ@>>>\dots @>>>C^nRZ@>>>\dots, \endCD $$
inducing the  isomorphism of the corresponding spectral sequences.  So
we have

{\bf Theorem 2.} {\sl The differentials of the Bousfield-Kan  spectral
sequence of the cosimplicial object
$$F^*(C,RZ):\quad RZ@>\delta^0,\delta^1>> CRZ\to\dots\to C^{n-1}RZ
@>\delta^0,\dots,\delta^n>>C^nRZ\to\dots$$ are    expressed   by   the
functional homology operations
$$H_*(\delta,\dots,\delta)\colon \pi_*(C^nRZ)\to   \pi_*(C^{n+m}RZ).$$
These operations determine on the $E^1$-term a new  differential.  The
homology  of  the  corresponding  chain  complex  is isomorphic to the
$E^\infty$-term of this spectral sequence.}

\vskip .5cm
\centerline{2. $E_\infty$-algebras and $E_\infty$-coalgebras}
\vskip 6 pt
Recall that  an  operad in the category of chain complexes is a family
$E=\{E(j)\}_{j\ge 1}$ of chain complexes $E(j)$  together  with  given
actions  of  symmetric groups $\Sigma_j$ and operations $$\gamma\colon
E(k)\otimes  E(j_1)\otimes\dots\otimes  E(j_k)\to   E(j_1+\dots+j_k)$$
compatible  with  these actions and satisfying associativity relations
[3], [4].

An operad  $E=\{E(j)\}$  for  which  complexes  $E(j)$ are acyclic and
symmetric groups act on them freely is called an $E_\infty$-operad.

A chain  complex $X$ is called an algebra (a coalgebra) over an operad
$E$ or simply an $E$-algebra (an $E$-coalgebra)  if  there  are  given
mappings  $$\mu\colon  E(k)\otimes_{\Sigma_k}X^{\otimes  k}\to  X\quad
(\tau\colon  X\to  Hom_{\Sigma_k}(E(k);X^{\otimes  k})),$$  satisfying
some associativity relations.

Algebras (coalgebras)    over    an   $E_\infty$-operad   are   called
$E_\infty$-algebras ($E_\infty$-coalgeb\-ras).

Any operad in the category  of  chain  complexes  determines  a  monad
$\underline E$ and a comonad $\overline E$ by the formulas

$$\gather \underline   E(X)=\sum_j\underline  E(j,X),\quad  \underline
E(j,X)=E(j)\otimes_{\Sigma_j}X^{\otimes        j};\\         \overline
E(X)=\prod_j\overline              E(j,X),\quad              \overline
E(j,X)=Hom_{\Sigma_j}(E(j),X^{\otimes j}).\endgather$$

An operad structure   $\gamma$ induces natural transformations
$$\underline\gamma\colon\underline  E\circ\underline E\to\underline E,
\qquad  \overline\gamma\colon\overline  E\to\overline  E\circ\overline
E.$$

An $E$-algebra  (an  $E$-coalgebra)  structure  on a chain complex $X$
induces a mapping $$\mu\colon \underline  E(X)\to  X\quad  (\tau\colon
X\to  \overline  E(X)).$$  So  to  give  on  a  chain  complex  $X$ an
$E$-algebra ($E$-coalgebra) structure is the same as to give on $X$ an
algebra  (coalgebra)  structure  over  the  monad  $\underline E$ (the
comonad $\overline E$).

One of the most important example  of  an  $E_\infty$-algebra  is  the
singular cochain complex $C^*(Y;R)$ of a topological space $Y$.

Dually, the singular chain complex $C_*(Y;R)$ of a topological space
$Y$ and the chain complex $N(RZ)$ of a simplicial set $Z$ are examples
of $E_\infty$-coalgebras.

The homotopy  theory  of $E_\infty$-coalgebras was constructed in [5].
There were defined the homotopy groups of  $E_\infty$-coalgebras.  For
the  chain  complex  $N(RZ)$  of  a  simplicial set $Z$ these homotopy
groups are isomorphic to the homotopy groups of $Z$ with  coefficients
in $R$.

For an   $E$-coalgebra   $X$,   using  cosimplicial  resolution
$$F^*(\overline  E,\overline  E,X):\quad  X@>\tau>>\overline   E(X)\to
\dots\to \overline E^{n-1}(X)\to\overline E^n(X)\to\dots, $$ there was
constructed the spectral  sequence  of  the  homotopy  groups  of  the
$E$-coalgebra $X$, [5].

Denote by $X_*$ the homology of the complex $X$ and by $\overline E_*$
the homology of a comonad $\overline  E$.  $X_*$  will  be  $\overline
E_*$-coalgebra.  There  is  the  cosimplicial resolution
$$F^*(\overline E_*,\overline E_*,X_*):\overline E_*(X_*)\to\overline
E_*^2(X_*)\to\dots\to\overline                 E_*^n(X_*)\to\overline
E_*^{n+1}(X_*)\to \dots.$$

The $E^1$ term of the spectral sequence  is  expressed  by  the  cobar
construction      $$F(\overline     E_*,X_*):\quad     X_*\to\overline
E_*(X_*)\to\dots\to          \overline          E_*^n(X_*)\to\overline
E_*^{n+1}(X_*)\to\dots,$$  obtained  from  the  resolution  by  taking
primitive elements. So there is the inclusion $F(\overline E_*,X_*)\to
F(\overline E_*,\overline E_*,X_*)$.

The functional  homology  operations  $$H_*(\delta,\dots,\delta)\colon
\overline  E^n_*(X_*)\to  \overline  E^{n+m}_*(X_*)$$  determine   new
differential  in  the  resolution  and in the cobar construction.  The
corresponding   complexes   denote    by    $\widetilde    F(\overline
E_*,\overline E_*,X_*)$, $\widetilde F(\overline E_*,X_*)$.

Note that  the complex $\widetilde F(\overline E_*,\overline E_*,X_*)$
is a resolution of the  complex  $X_*$  and  there  is  the  inclusion
$\widetilde    F(\overline    E_*,X_*)\to    \widetilde    F(\overline
E_*,\overline E_*,X_*)$.

{\bf Theorem 3.} {\sl Differentials of the spectral  sequence  of  the
homotopy groups   of  an  $E$-coalgebra  $X$  are  determined  by  the
functional homology operations
$$H_*(\delta,\dots,\delta)\colon  \overline  E^n_*(X_*)\to   \overline
E^{n+m}_*(X_*).$$ The homology of $\widetilde F(\overline E_*,X_*)$ is
isomorphic to the $E^\infty$-term of the spectral sequence.}

If $X$  is the normalized chain complex of a simplicial set $Z$,  i.e.
$X=N(RZ)$, then there is a mapping of cosimplicial objects
$$\CD
N(RZ)@>>>N(CRZ)@>>>\dots  @>>>N(C^nRZ)@>>>\dots\\  @V=VV @VVV @.  @VVV
@.\\ X@>>>\overline E(X)@>>>\dots @>>>\overline E^n(X)@>>>\dots,\endCD
$$ inducing the isomorphism of the corresponding spectral sequences.
So we have

{\bf Theorem 4.} {\sl The differentials of the Bousfield-Kan spectral
sequence of the homotopy groups of a simplicial set $Z$ are determined
by the functional homology operations $$H_*(\delta,\dots,\delta)\colon
\overline E^n_*(Z_*)\to \overline E^{n+m}_*(Z_*).$$
The homology of $\widetilde F(\overline E_*,Z_*)$ is isomorphic to the
$E^\infty$-term of the spectral sequence.}

Note that  the  suspension  $SX$  over  an  $E$-coalgebra  $X$  is  an
$SE$-coalgebra and the following diagrams commute
$$\CD \overline  E@>\overline\gamma  >>\overline  E\circ\overline  E\\
@VVV                                                            @VVV\\
\overline{SE}@>\overline{S\gamma}>>\overline{SE}\circ\overline{SE}
\endCD\qquad  \CD  SX@>\tau  >>   \overline{SE}(SX)\\   @V=VV   @VVV\\
SX@>S\tau >>S(\overline E(X))\endCD $$ Moreover from the expression of
the homology $\overline E_*$ of the comonad $\overline E$ (see  below)
it  follows  that the mappings $\xi\colon\overline E_*\to\overline E$,
$\eta\colon\overline  E\to\overline   E_*$,   $h\colon\overline   E\to
\overline E$ may be chosen permutable with the suspension homomorphism
$\overline{SE}\to\overline  E$.   Therefore   constructed   functional
homology  operations  permute with the suspension homomorphism.  Hence
the following theorem is taken place.

{\bf Theorem  5.}  {\sl  Functional homology operations giving higher
differentials of the Bousfield-Kan spectral sequence permute with  the
suspension and hence are stable.  They induce the differentials of the
Adams spectral sequence of stable homotopy  groups  of  a  topological
space.}

\vskip .5cm
\centerline{3. $A_\infty$-cosimplicial objects}
\vskip 6pt
Make more   precise   the   form   of   higher  differentials  of  the
Bousfield-Kan spectral   sequence   using    the    notion    of    an
$A_\infty$-cosimplicial object.

A family $X^*=\{X^n\}_{n\ge 0}$ of objects $X^n$ of a category $\Cal X$
will be called a precosimplicial object if there are given coface  and
codegeneracy operators
$$\gather\delta^i\colon X^n\to X^{n+1},\quad 1\le  i\le
n+1;\\  \sigma^i\colon X^n\to X^{n-1},\quad 0\le i\le n-1,\endgather$$
satisfying the following relations
$$\align
\delta^j\delta^i&=\delta^i\delta^{j-1},~                        i<j,\\
\sigma^j\sigma^i&=\sigma^i\sigma^{j+1},~i\le                      j,\\
\sigma^j\delta^i&=\cases                  \delta^i\sigma^{j-1},&i<j,\\
Id,&i=j,~i=j+1,\\ \delta^{i-1}\sigma^j,&i>j+1.\endcases \endalign $$

Thus a precosimplicial object differs a cosimplicial objects only by a
coface operators $\delta^0$.
A cosimplicial   object   has   such  operator  but  a  precosimlicial
object hasn't.

A mapping $f^*\colon X\to Y$ of a precosimplicial objects is a family
$f^*=\{f^n\}_{n\ge 0}$  of  mappings  $f^n\colon X^n\to Y^n$ commuting
with coface and codegeneracy operators,  i.e. satisfying the following
relations
$$\align
\delta^if_n&=f^{n+1}\delta^i,\\ \sigma^if^n&=f_{n-1}\sigma^i.\endalign
$$

Now we define the notion of an $A_\infty$-cosimplicial object  of  the
category of topological spaces.

Let $I^n$ be the unite cube, $I^n=\{(t_1,\dots,t_n)|~0\le t_i\le 1\}$.
Denote by
$$\align
u_i^\epsilon\colon I^n\to I^{n+1},&~\epsilon =0,1;~1\le i\le n+1;\\
v_i\colon I^n\to I^{n-1},&~ 0\le i\le n,\endalign $$
the mappings defined by the formulas
$$\align u_i^\epsilon(t_1,\dots,t_n)&=
(t_1,\dots t_{i-1},\epsilon,t_i,\dots,t_n);\\
v_i(t_1,\dots,t_n)&=\cases (t_2,\dots,t_n),&i=0,\\
(t_1,\dots,t_i*t_{i+1},\dots,t_n),&1\le i\le n-1,\\
(t_1,\dots,t_{n-1}),& i=n,\endcases \endalign $$
where $t_i*t_{i+1}=t_i+t_{i+1}-t_i\cdot t_{i+1}$.

A precosimplicial object $X^*=\{X^n\}$ of the category of  topological
spaces will  be  called  an  $A_\infty$-cosimplicial  object or simply
an $A_\infty$-cosimplicial space if there are given coface operators
$$\delta^0_m\colon  X^n\times I^m\to X^{n+m+1},$$
satisfying the relations
$$\align
\sigma^0\delta^0_0&=Id;\\                           \delta^0_m(1\times
u_i^0)&=\delta^i\delta^0_{m-1},~1\le   i\le   m;\\  \delta^0_m(1\times
u_i^1)&=\delta^0_{i-1}(\delta^0_{m-1}\times  1),   ~1\le   i\le   m;\\
\delta^0_{m-1}(1\times   v_i)&=\sigma^i\delta^0_m,~0\le  i\le  m,~m\ge
1;\\ \delta^0_m(\delta^i\times 1)&=\delta^{i+m+1}\delta^0_m,~i\ge 1,\\
\delta^0_m(\sigma^i\times     1)&=\sigma^{i+m+1}\delta^0_m,~i\ge    0.
\endalign $$

It is clear that a cosimplicial object $X^*=\{X^n\}$ of  the  category
of topological spaces may be concidered as an $A_\infty$-cosimplicial
object with trivial operators
$\delta^0_m\colon X^n\times I^m\to X^{n+m+1}$ при $m\ge 1$.

Note that  the  family  $I^*=\{I^n\}$  itself may be considered as an
$A_\infty$-cosimplicial space for which
$$\gather\delta^i=u_i^0\colon
I^n\to  I^{n+1},\quad  1\le  i\le  n+1;\\   \delta^0_m=u_{n+1}^1\colon
I^n\times   I^m=I^{n+m}\to   I^{n+m+1};\\   \sigma^i=v_i\colon  I^n\to
I^{n-1},\quad 0\le i\le n.\endgather $$

Define also    the    notion    of    an     $A_\infty$-mapping.

Let $X^*=\{X^n\}$   be   an    $A_\infty$-cosimplicial    space    and
$Y^*=\{Y^n\}$ be  a  cosimplicial space.  Then $A_\infty$-mapping from
$X^*$ to $Y^*$  is a precosimplicial mapping
$f^*=\{f^n\}$, $f^n\colon  X^n\to  Y^n$  together  with  the family of
mappings
$$f^n_m\colon X^n\times    I^m\to   Y^{n+m},\quad   0\le   m\le   n,$$
satisfying the following relations
$$\align    f^n_0&=f^n;\\
f^n_m(1\times     u_i^0)&=\delta^{i-1}f^n_{m-1},~1\le     i\le    m;\\
f^n_m(1\times u_i^1)&=f^{n+m-i-1}_i(\delta^0_{m-i}\times 1),~1\le i\le
m;\\  f_{n-1}^m(d_i\times  1)&=d_{i-m}f_n^m,~i>m;\\  f^n_{m-1}(1\times
v_i)&=\sigma^if^n_m,~0\le         i<m;\\          f^n_m(\delta^i\times
1)&=\delta^{i+m}f^{n-1}_m,~i\ge      1;\\     f^n_{m-1}(\sigma^i\times
1)&=\sigma^{i+m}f_n^m,~i\ge 0. \endalign $$

An $A_\infty$-mapping    will   be   called   an   $A_\infty$-homotopy
equivalence if  the  corresponding   mappings   $f_n$   are   homotopy
equivalences.

Of course  the  notion  of  an  $A_\infty$-cosimplicial  object may be
defined not only for the category of topological spaces  but  also  in
some other situations.  For example it may be defined for the category
of simplicial sets,  for the category of chain complexes and so on. To
do it  we need to use the analog of the $n$-dimensional cube $I^n$ for
these categories.

Consider more  precisely  $A_\infty$-cosimplicial  objects   for   the
category of  chain  complexes.  In the capacity of $I^n$ we take the
normalized  chain  complex of $n$-dimensional cube.

The definitions   of    an    $A_\infty$-cosimplicial    object    and
an $A_\infty$-mappings may be reformulated in the following form.

A precosimplicial  object  $X^*=\{X^n\}$  of  the  category  of  chain
complexes is  an   $A_\infty$-cosimplicial   object   or   simply   an
$A_\infty$-cosimplicial complex   if  there  are  given  mappings
$$\delta^0_m\colon X^n\to X^{n+m+1},$$ increasing  dimensions  by  $m$
and satisfying the following relations $$\align
     d(\delta^0_m)=&\sum_{i=1}^n(-1)^{i-1}(\delta^i\delta^0_{m-1}-
\delta^0_{i-1}\delta^0_{m-i});\\
     \delta^0_m\delta^i=&\delta^{i+m+1}\delta^0_m,~i\ge1;\\
     \sigma^i\delta^0_m=&\cases 0,&0\le i\le m,~m\ge 1,\\
     Id,&i=m=0,\\\delta^0_m\sigma_{i-m-1},&i>m.\endcases\endalign $$

An $A_\infty$-mapping from an $A_\infty$-cosimplicial complex $X^*$  to
a cosimplicial    complex   $Y^*$   is   a   precosimplicial   mapping
$f^*=\{f^n\}$,  $f^n\colon X^n\to Y^n$ together  with  a family  of
mappings $f^n_m\colon X^n\to Y^{n+m}$ increasing dimensions by $m$ and
satisfying the  following  relations     $$\align       f^n_0=&f^n;\\
d(f_m^n)=&\sum_{i=1}^n(-1)^{i-1}(\delta^{i-1}f^n_{m-1}-f^{n+m-i-1}_i
\delta^0_{m-i});\\   f_m^{n+1}\delta^i=&\delta^{i+m}f^n_m,~i\ge   1;\\
\sigma^if_m^n=&\cases 0,&i<m;\\f^{n+1}_m\sigma^{i-m},&i\ge m.\endcases
\endalign $$

It is  clear  that  a  cosimplicial  complex  $X^*=\{X^n\}$   may   be
considered as   an  $A_\infty$-cosimplicial   complex  with  trivial
operators $\delta^0_m\colon X^n\otimes I^m\to X^{n+m+1}$ if $m\ge 1$.

Let ${X'}^*$,  ${X''}^*$ are $A_\infty$-cosimplicial  complexes.  Then
there there is  defined  the tensor product ${X'}^*\otimes {X''}^*$.
It is an $A_\infty$-cosimplicial complex $X^*$ for which
$$\gather
X^n={X'}^n\otimes                                           {X''}^n;\\
\delta^i={\delta'}^i\otimes{\delta''}^i,\quad        i\ge         1;\\
\sigma^i={\sigma'}^i\otimes{\sigma''}^i,\quad  i\ge  0.\endgather $$
Operations $\delta^0_m\colon  X^n\to  X^{n+m+1}$ are the compositions
$$\gather           {X'}^n\otimes{X''}^n\otimes          I^m@>1\otimes
1\otimes\nabla>>{X'}^n \otimes{X''}^n\otimes I^m\otimes I^m@>T>>\\ \to
{X'}^n\otimes    I^m\otimes{X''}^n\otimes    I^m@>{\delta'}^0_m\otimes
{\delta''}^0_m>>{X'}^{n+m+1}\otimes{X''}^{n+m+1},\endgather   $$ where
$\nabla\colon I^m\to I^m\otimes I^m$ is a comultiplication  in  the
coalgebra $I^m$.

Given $A_\infty$-cosimplicial  complex  $X^*$  define  its realization
$|X^*|$ putting $$|X^*|=Hom(I^*;X^*),$$ where $Hom$ is  considered  in
the category of $A_\infty$-cosimplicial complexes.

Similary as    for     the     usual     realization     for     given
$A_\infty$-cosimplicial complexes ${X'}^*$, ${X''}^*$ there is a chain
equivalence    $$|{X'}^*\otimes     {X''}^*|\simeq     |{X'}^*|\otimes
|{X''}^*|.$$

Transfer the  Perturbation  Theory  [11]  to  $A_\infty$-cosimplicial
objects.

We say that a chain complex $\widetilde X$ is a deformation  retract
of a chain complex $X$ if there are chain mappings
$\xi\colon\widetilde X\to X$,  $\eta\colon  X\to\widetilde  X$  and  a
chain homotopy $h\colon  X\to  X$ such that
$$ \eta\circ\xi=Id,\quad d(h)=Id-\xi\circ\eta.$$ In this case  we  may
additionally assume that
$$h\circ\xi=0,\quad\eta\circ h=0,\quad h\circ h=0,\quad [10].$$

A precosimplicial  complex  $\widetilde  X^*$   will   be   called   a
deformation retract  of a precosimplicial complex $X^*$ if for any $n$
a chain complex $\widetilde X^n$ is a deformation retract of a complex
$X^n$ and  the  corresponding  chain  mappings  and  chain  homotopies
commute with a precosimplicial structure.

The next theorem is analog of  the  main  lemma  of  the  Perturbation
Theory [11].

{\bf Theorem  6.} {\sl Let $X^*=\{X^n\}$ be a cosimplicial complex and
a precosimplicial complex $\widetilde  X^*=\{\widetilde  X^n\}$  is  a
deformation retract  of $X^*$ considered as a precosimplicial complex.
Then on $\widetilde X^*$ there  is an $A_\infty$-cosimplicial  structure
and an $A_\infty$-cosimplicial homotopy equivalence between $X^*$ and
$\widetilde X^*$.}

{\bf Proof.}  Let  $\xi^*\colon\widetilde   X^*\to  X^*$,
$\eta^*\colon X^*\to \widetilde X^*$ and $h^*\colon X^*\to X^*$ be the
corresponding mappings.  Define operators $\delta^0_m\colon \widetilde
X^n\to\widetilde                  X^{n+m+1}$ putting
$$\delta^0_m=\eta\delta^0h^{n+m}\dots    \delta^0h^{n+1}\delta_0\xi.$$
Direct calculations show that the required relations are satisfied.

Define also mappings $\xi^n_m\colon \widetilde X^n\to X^{n+m}$ putting
$$\xi^n_m=h^{n+m}\delta^0\dots      h^{n+1}\delta^0\xi.$$       Direct
calculations    show    that    these    mappings    give    us    the
$A_\infty$-cosimplicial homotopy equivalence between $\widetilde  X^*$
and $X^*$.

{\bf Theorem  7.}  {\sl  Let  $f^*\colon  X^*\to  Y^*$ be a mapping of
cosimplicial complexes and precosimplicial complexes $\widetilde X^*$,
$\widetilde Y^*$ are deformation retracts of $X^*$,  $Y^*$, considered
as    precosimplicial    complexes.    Then    for    above    defined
$A_\infty$-cosimplicial  structures  on $\widetilde X^*$,  $\widetilde
Y^*$ there is an  $A_\infty$-mapping  $\widetilde  f^*\colon\widetilde
X^*\to\widetilde  Y^*$  such that the following diagram commutes up to
homotopy
     $$\CD X^*@>f^*>> Y^*\\
     @A\xi^*AA @AA\xi^*A\\
     \widetilde X^*@>\widetilde f^*>>\widetilde Y^*\endCD$$}

The proof is analog to the previous theorem.

Taking into account that the $E^1$-term of the Bousfield-Kan  spectral
sequence, considered  as  a  precosimplicial object, is a deformation
retract of the initial cosimplicial object we obtain

{\bf Theorem 8.} {\sl The functional homology operations giving higher
differentials  of the Bousfield-Kan spectral sequence may be chosen in
such a way  to  form  on  the  $E^1$-term  an  $A_\infty$-cosimplicial
structure.}

\vskip .5cm
\centerline{4. $E_\infty$-structure on the Bousfield-Kan spectral
sequence}
\vskip 6pt
Our aim here is to define $E_\infty$-structure  on  the  Bousfield-Kan
spectral sequence.  To do it consider the following additional propety
of the operad $E$.

{\bf Theorem  9.}  {\sl  Given  $E_\infty$-operad  $E$  there   is   a
permutation  mapping $$T\colon\underline E\circ\overline E\to\overline
E\circ  \underline  E,$$  commuting  the  following   diagrams   $$\CD
\underline   E^2\circ\overline   E@>T\underline  E\circ\overline  ET>>
\overline  E\circ\underline  E^2\\  @V\underline\gamma\overline  E  VV
@VV\overline    E\underline\gamma   V\\   \underline   E\circ\overline
E@>T>>\overline  E\circ\underline   E   \endCD   \quad\CD   \underline
E\circ\overline   E@>T>>\overline  E\circ\underline  E\\  @V\underline
E\overline\gamma  VV  @VV\overline\gamma\underline  E  V\\  \underline
E\circ\overline   E^2@>\overline   ET\circ  T\overline  E>>  \overline
E^2\circ\underline E\\ \endCD $$}

{\bf Proof.} Given an $E_\infty$-operad $E$ we may construct an operad
mapping  $\nabla\colon  E\to  E\otimes  E$  consisting   of   mappings
$\nabla(j)\colon  E(j)\to E(j)\otimes E(j)$ This mapping give on $E$ a
Hopf operad structure. Denote $\nabla(j,i)\colon E(j)\to E(j)^{\otimes
i}$ the iterations of these mappings $\nabla(j,2)=\nabla(j)$. They are
$\Sigma_j$-mappings,                                              i.e.
$$\nabla(j)(x\sigma)=\nabla(j)(x)\sigma^{\otimes     j},\quad\sigma\in
\Sigma_j$$ but are not  commuting  with  permutations  of  factors  of
$E(j)^{\otimes  i}$.  However  they  may be extended till the mappings
$\nabla(j,i)\colon E(i)\otimes E(j)\to E(j)^{\otimes  i}$,  compatible
with the actions of symmetric groups $\Sigma_i$ and $\Sigma_j$.

Rewrite these  mappings  in  the  form   $$\nabla(j,i)\colon   E(j)\to
\overline  E(i)\otimes E(j)^{\otimes i}.$$ If the operad $E$ is chosen
freely then they may be done compatible with an operad structure.

Passing to    the     dual     mappings     we     obtain     mappings
$$\overline\nabla(j,i)\colon  E(i)\otimes\overline E(j)^{\otimes i}\to
\overline E(j).$$

Define mappings   $$T(j,i)\colon   E(j)\otimes\overline  E(i)^{\otimes
j}\to E(i)\otimes\overline E(j)^{\otimes i}$$ as the compositions
     $$\gather
E(j)\otimes\overline E(i)^{\otimes   j}@>\nabla(j)\otimes   1^{\otimes
j}>>E(j)\otimes E(j)\otimes\overline  E(i)^{\otimes   j}\to   \\   \to
E(j)\otimes\overline              E(i)^{\otimes              j}\otimes
E(j)@>\overline\nabla(j,i)\otimes 1>>\overline E(i)\otimes E(j)\to  \\
@>1\otimes\nabla(j,i)>>\overline E(i)\otimes\overline      E(i)\otimes
E(j)^{\otimes i} @>\overline\nabla(i)\otimes 1^{\otimes  i}>>\overline
E(i)\otimes E(j)^{\otimes i}. \endgather $$

     The family  $T(j,i)$ determines the required permutation mapping
$$T\colon\underline E\circ\overline E\to\overline E\circ\underline E.$$

A chain complex $X$ will be called an $E_\infty$-Hopf algebra if there
are given an $E_\infty$-algebra structure $\mu\colon\underline E(X)\to
X$ and an  $E_\infty$-coalgebra  structure  $\tau\colon  X\to\overline
E(X)$  such  that  the  following  diagram  commutes  $$\CD \underline
E(X)@>\mu>>X@>\tau>>\overline E(X)\\ @V\underline E(\tau)VV @. @VV=V\\
\underline  E\overline E(X)@>T>>\overline E\underline E(X) @>\overline
E(\mu)>>\overline E(X)\endCD $$

If a  topological  space  $Y$ is an $E_\infty$-space then its singular
chain complex $C_*(Y;R)$  will  be  an  $E_\infty$-Hopf  algebra.  For
example the singular chain complex of the infinite loop space is an
$E_\infty$-Hopf algebra.

Let $X$  is  an  $E_\infty$-Hopf  algebra.  Then there is a mapping of
augmented cosimplicial  objects  $$\CD  \underline  E(X)@>>>\underline
E\overline  E(X)@>>>\dots @>>> \underline E\overline E^n(X)@>>>\dots\\
@VVV @VVV @.  @VVV  @.\\  X@>>>\overline  E(X)@>>>\dots  @>>>\overline
E^n(X)@>>>\dots\endCD $$

Then the   complex   $F(\overline   E,\overline    E,X)$    will    be
$E_\infty$-algebra. Passing     to     the    homology    we    obtain
$E_\infty$-algebra structure on the  complex  $\widetilde  F(\overline
E_*,\overline E_*,X_*)$. Thus we have

{\bf Theorem 10.} {\sl If $X$ is an $E_\infty$-Hopf algebra,  then the
complex $\widetilde  F(\overline  E_*,\overline  E_*,X_*)$   possesses
$E_\infty$-algebra structure.}

This structure  will  be  used  in  the further calculations of higher
differentials of the Bousfield-Kan spectral sequence.

Note that on the cobar construction $\widetilde F(\overline  E_*,X_*)$
there is no $E_\infty$-algebra structure.

\vskip .5cm
\centerline{5. The  homology  of  an  $E_\infty$-operad and the Milnor
coalgebra}
\vskip 6pt
Let $E$ be an $E_\infty$-operad,  $M$ be a graded module  (over  $\Bbb
Z/2$). As it is known (see for example [6]), the homology $\underline
E_*(M)$ of the complex $\underline E(M)$  is  the  polinomial  algebra
generated by    the    elements   $e_{i_1}\dots   e_{i_k}x_m$,   $1\le
i_1\le\dots\le i_k$,  $x_m\in   M$   of   dimensions  $i_1+2i_2+\dots+
2^{k-1}i_k+2^km$.

The elements  $e_{i_1}\dots  e_{i_k}x_m$ of $\underline E_*(M)$ may be
rewritten in the form
$$Q^{j_1}\dots Q^{j_k}\otimes  x_m;\quad  j_1\le 2j_2,\dots,j_{k-1}\le
2j_k,          m\le          j_k,$$
where
$$\gather j_k=i_k+m,\\j_{k-1}=i_{k-1}+i_k+2m,\\             \dots             \\
j_1=i_1+i_2+2i_3+\dots+2^{k-2}i_k+2^{k-1}m.\endgather $$
The sequences  $Q^{j_1}\dots  Q^{j_k}$  give  the  elements   of   the
Dyer-Lashof algebra $\Cal R$, [12], [13].

Given a  graded  module  $M$  denote  by $\Cal R\times M$ the quotient
module of the tensor product $\Cal R\otimes  M$  under  the  submodule
generated by the elements $Q^{j_1}\dots Q^{j_k}\otimes x_m$,  $j_k<m$.
The correspondence $M\longmapsto \Cal R\times M$  determines the monad
in the category of graded modules.

A graded  module $M$ is called an unstable module over the Dyer-Lashof
algebra if it is an algebra over the corresponding monad.

Dually, the homology $\overline  E_*(M)$  of  the  complex  $\overline
E(M)$ is the free commutative coalgebra generated by the elements
$$e_{i_1}\dots e_{i_k}x^m,\quad 1\le i_1\le\dots\le  i_k,\quad  x^m\in
M$$ of dimensions $2^km-(i_1+2i_2+\dots+ 2^{k-1}i_k)$.

Regrading of  the  elements  of $\overline E_*(M)$ leads to the Milnor
coalgebra $\Cal  K$.
By definition $\Cal K$ is the  polynomial  algebra  generated  by  the
elements $\xi_i$, $i\ge 0$ of dimensions $2^i-1$. A   comultiplication
$$\nabla\colon  \Cal  K\to\Cal  K\otimes\Cal  K$$  on  the  generators
$\xi_i$         is         given         by         the        formula
$$\nabla(\xi_i)=\sum_k\xi_{i-k}^{2^k}\otimes\xi_k.$$  On   the   other
elements the comultiplication is determined by the Hopf relation.

Define the   grading   $deg(x)$   of  elements  $x\in\Cal  K$  putting
$deg(\xi_i)=1$ and the grading of the product equal to the sum of  the
gradings of factors.

Given graded  module  $M$  denote by $\Cal K\times M$ the submodule of
the tensor  product  $\Cal  K\otimes  M$  generated  by  the  elements
$x\otimes y$,  $deg(x)=dim(y)$.  The correspondence $M\longmapsto \Cal
K\times M$ determines the comonad $\Cal K$ in the category  of  graded
modules.

A graded module $M$ is called an unstable  comodule  over  the  Milnor
coalgebra if it is a coalgebra over the corresponding comonad.

Let $M$ be an unstable module over the Milnor coalgebra.  There is a
cosimplicial resolution $$F^*(\Cal K,\Cal K,M):M\to\Cal K\times
M\to\dots\to\Cal  K^{\times  n}  \times M \to\Cal K^{\times n+1}\times
M\to\dots $$ If $Y$ is a "nice"  (in  the  sense  of  Massey-Peterson)
space then   the   Bousfield-Kan   spectral   sequence  turns  to  the
Massey-Peterson spectral sequence.  The $E^1$-term  of  this  spectral
sequence may be written in the form
$$F^*(\Cal K,Y_*):Y_*\to\Cal  K\times  Y_*\to\dots\to  \Cal  K^{\times
n}\times Y_* \to\dots, $$ where $Y_*=H_*(Y;\Bbb Z/2)$.

From the previos theorems it follows the next theorem.

{\bf Theorem  11.} {\sl The functional homology operations determining
higher differentials of the Massey-Peterson  spectral  sequence  of  a
"nice"  space  $Y$  may  be  chosen in such a way to form on $F^*(\Cal
K,\Cal K,Y_*)$ an $A_\infty$-cosimplicial structure.  The homology  of
the  corresponding  cobar  construction  $\widetilde F(\Cal K,Y_*)$ is
isomorphic to the $E^\infty$  term  of  the  Massey-Peterson  spectral
sequence.  If  $Y$ is an $E_\infty$-space then the complex $\widetilde
F(\Cal K,\Cal K,Y_*)$ is $E_\infty$-algebra.}

Note that  there  is  an  inclusion   $\widetilde   F(\Cal   K,Y_*)\to
\widetilde  F(\Cal  K,\Cal K,Y_*)$.  However on the $\widetilde F(\Cal
K,Y_*)$ there is no $E_\infty$-algebra structure.

Besides the Milnor coalgebra $\Cal K$  we  will  consider  the  stable
Milnor coalgebra  $\Cal  K_s$ for which $\xi_0=1$.

Given comodule  $M$  over  the  stable  Milnor  coalgebra  there is a
cosimplicial resolution
$$F^*(\Cal  K_s,\Cal K_s,M):\Cal K_s\otimes M\to\Cal  K_s^{\otimes 2}
\otimes M\to\dots\to\Cal K_s^{\otimes n} \otimes M\to\dots $$

Stabilization of the Bousfild-Kan spectral sequence leads to the Adams
spectral sequence of stable homotopy groups of a topological space $Y$.
$E^1$-term of this spectral sequence may be written in the form
$$F^*(\Cal K_s,Y_*):Y_*\to\Cal    K_s\otimes    Y_*\to\dots\to    \Cal
K_s^{\otimes n}\otimes Y_*\to\dots $$ Thus we have

{\bf Theorem  12.}  {\sl  Functional  homology  operations determining
higher differentials of the Adams spectral sequence of stable homotopy
groups  of a topological space $Y$ may be chosen in such a way to form
on   $F^*(\Cal   K_s,\Cal    K_s,Y_*)$    the    structure    of    an
$A_\infty$-cosimplicial  object.  The  homology  of  the corresponding
cobar construction $\widetilde F(\Cal K_s,Y_*)$ is isomorphic  to  the
$E^\infty$ term of the Adams spectral sequence.}

Let us  calculate  the  $E_\infty$-algebra  structure  on  the  Milnor
coalgebra.  As it was pointed out above for an  $E_\infty$-operad  $E$
there  is  the  permuting  mapping  $T\colon\underline E\circ\overline
E\to\overline E\circ \underline E$.  It induces the permuting  mapping
$T_*\colon\underline    E_*\circ\overline    E_*\to    \overline   E_*
\circ\underline   E_*$    commuting    diagrams    $$\CD    \underline
E^2_*\circ\overline    E_*@>T_*\underline   E_*\circ\overline   E_*T_*
>>\overline E_*\circ\underline  E^2_*\\  @V\underline\gamma_*\overline
E_*VV       @VV\overline      E_*\underline\gamma_*V\\      \underline
E_*\circ\overline E_*@>T_*>>  \overline  E_*\circ\underline  E_*\endCD
\qquad    \CD    \underline    E_*\circ\overline   E_*@>T_*>>\overline
E_*\circ\underline    E_*\\    @V\underline     E_*\overline\gamma_*VV
@VV\overline\gamma_*  \underline  E_*V\\  \underline E_*\circ\overline
E_*^2@>\overline     E_*T_*\circ      T_*\overline      E_*>>\overline
E_*^2\circ\underline E_* \endCD $$

The permuting mapping $T_*$ induces the action $\mu_*\colon \underline
E_*\circ  \overline  E_*\to  \overline  E_*$  and  the  dual  coaction
$\tau_*\colon \underline E_*\to\overline E_*\circ\underline E_*$.

Denote by  $e_i\colon\Cal  K\to\Cal  K$  the  operation  on the Milnor
coalgebra inducing by the restriction of $\mu_*$ on the elements $e_i$.
From the  commutative  diagrams  for  the  permuting  mapping $T_*$ it
follows

{\bf Theorem  13.}  {\sl The operations  $e_i\colon\Cal  K\to\Cal  K$
satisfy the relations:

1. $e_0(x) = x^2$.

2. $e_i(xy) = \sum e_k(x)e_{i-k}(y)$.

3. $\nabla e_i(x) = \sum \xi_0^{-k}e_{i-k}(\xi_0^{k}x')
\otimes e_k(x'')$, where $\sum x'\otimes x''=\nabla(x) $.}

     Using these  relations  to  calculate  the operations $e_i$ it is
sufficient to calculate only $e_1(\xi_0)$.  Direct  calculasions  show
that $e_1(\xi_0)=\xi_1\xi_0$. From the third relation it follows

{\bf Theorem 14.}  {\sl  There are the following formulas
$$  e_i(\xi_k)  =  \cases
\xi_{m+k}\xi_k,   &   i  =2^{m+k}-2^k;\\  \xi_{m+k}\xi_{k-1},  &  i  =
2^{m+k}-2^k-2^{k-1};\\  \dots  &  \dots  \\   \xi_{m+k}\xi_0,&   i   =
2^{m+k}-2^k-\dots -1;\\ 0,& \text{in other cases.} \endcases $$}

     Using the  second  relation  we  may  obtain  formulas  for   the
operations $e_i$ on the products of the elements $\xi_k$.

Passing from the elements $e_i$ to the elements  of  the  Dyer-Lashof
algebra we  obtain the action of the Dyer-Lashof algebra on the Milnor
coalgebra. On the generators  $\xi_i$ it is given by the formulas
$$  Q^{i+2^k-1}(\xi_k)  =
\cases \xi_{m+k}\xi_k,  & i =2^{m+k}-2^k;\\ \xi_{m+k}\xi_{k-1},  & i =
2^{m+k}-2^k-2^{k-1};\\  \dots  &  \dots  \\   \xi_{m+k}\xi_0,&   i   =
2^{m+k}-2^k-\dots -1;\\ 0,& \text{in other cases.} \endcases $$

On the other elements this action is determined by the Hopf relations
$$Q^i(xy)=\sum Q^k(x)Q^{i-k}(y).$$

Besides the  action of the Dyer-Lashof algebra on the Milnor coalgebra
$\Cal K$  there  are  $\cup_i$-products  and   an   $E_\infty$-algebra
structure. On  the  generators  $x\in\Cal  K$  it  is  defined  by the
formulas $$x\cup_ix=e_i(x).$$ On the other elements  $\cup_i$-products
are defined by the relations
$$\gather x\cup_i                                        y=y\cup_ix,\\
(x_1x_2)\cup_iy=x_1(x_2\cup_iy)+(x_1\cup_iy)x_2.\endgather$$

Note that  the stable Milnor coalgebra $\Cal K_s$ has no action of the
Dyer-Lashof algebra and has no $E_\infty$-algebra structure.

\vskip .5cm
\centerline{6. Degenerating $A_\infty$-structures}
\vskip 6pt
As it  was  pointed  out  above  the  Milnor  coalgebra  $\Cal  K$  is
$E_\infty$-algebra and in particular $A_\infty$-algebra.  Here we show
that this $A_\infty$-atructure on the Milnor coalgebra is degenerated.

Recall that a chain complex $A$ is called  an  $A_\infty$-algebra  [2]
(over   $\Bbb  Z/2$)  if  there  are  given  operations  $$\pi_n\colon
A^{\otimes n+2} \to A,\quad n\ge 0,$$ increasing dimensions by $n$ and
satisfying      the      following      relations      $$d(\pi_{n+1})=
\sum_{i=0}^n\pi_i(1\otimes\dots\otimes\pi_{n-i}\otimes\dots    \otimes
1),$$ where the sum is taken also overe all places of $\pi_{n-i}$.

In particular a graded module $A_*$ is an $A_\infty$-алгеброй if there
are given operations
$$\pi_n\colon A_*^{\otimes n+2} \to A_*,\quad n\ge 0,$$
increasing dimensions by $n$ and satisfying  the  following  ralations
$$\sum_{i=0}^n\pi_i(1\otimes\dots\otimes\pi_{n-i}\otimes\dots\otimes
1)=0.$$

Let $A'$,  $A''$  be  $A_\infty$-algebras.  A   family   of   mappings
$$f_n\colon  {A'}^{\otimes  n+1}\to  A'',\quad  n\ge  0,$$  increasing
dimensions by $n$ is called an $A_\infty$-mapping from $A'$  to  $A''$
if     the     following     relations    are    satisfied    $$\align
d(f_{n+1})&=\sum_{i=0}^nf_i(1\otimes\dots\otimes\pi'_{n-i}\otimes
\dots\otimes                                                     1)-\\
&-\sum_{n_1+\dots+n_{i+2}=n-i}\pi''_i(f_{n_1}\otimes\dots\otimes
f_{n_{i+2}}).\endalign$$

Two $A_\infty$-structures $\{\pi'_n\}$ and  $\{\pi''_n\}$  are  called
equivalent    if   there   exists   an   $A_\infty$-mapping   $f\colon
(A,\pi')\to(A,\pi'')$, for which $f_0=Id\colon A\to A$.

An $A_\infty$-structure is called degenerated if it is  equivalent  to
usual algebra structure.

$A_\infty$-algebra structure appears by a natural way on the  homology
$A_*$  of  a  differential  algebra  $A$,  [14].  Namely,  we  fix the
homomorphism $\xi\colon  A_*\to  A$  of  choosing  representatives  in
homology  classes,  fix the inverse homomorphism $\eta\colon A\to A_*$
and   a    chain    homotopy    $h\colon    A\to    A$    such    that
$$\eta\circ\xi=Id,\quad d(h)=\xi\circ\eta-Id,  \quad h\circ\xi=0,\quad
\eta\circ h=0,\quad h\circ h=0.$$ Then $A_\infty$-algebra structure on
$A_*$  is  given by the formula $$\pi_n=\eta\pi(h\pi\otimes 1+1\otimes
h\pi)\dots(h\pi\otimes\dots\otimes    1    +\dots+1\otimes\dots\otimes
h\pi)(\xi\otimes\dots\otimes\xi).$$

It is easy to see that if the homology $A_*$ of a differential algebra
$A$ is isomorphic to the tensor algebra $TX$  generated  by  a  graded
module $X$, then $A_\infty$-algebra structure on $A_*$ is degenerated.
In this case there is an algebra mapping $\xi\colon A_*\to  A$  giving
by          the          formula          $$\xi(x_1\otimes\dots\otimes
x_n)=\xi(x_1)\cdot\dots\cdot \xi(x_n),$$ where $x_i\in X$,  $\xi(x_i)$
- representatives of elements $x_i$.

Consider the  case when $A_*$ is the polynomial algebra $PX$ generated
by a graded module $X$.

If $A$ is a  commutative  algebra  then  $A_\infty$-algebra  structure
on $A_*=PX$ degenerated. The corresponding algebra mapping
$\xi\colon  A_*\to  A$ is given by the formula $$\xi(x_1\cdot\dots\cdot
x_n)=\xi(x_1)\cdot\dots\cdot \xi(x_n).$$

In the  case  of a none commutative algebra $A$ the mapping $\xi\colon
A_*=PX\to A$ may be choosen in such a way that $\xi(x\cdot y)=\xi(x)
\cdot\xi(y)$ if $x\le y$.

From the  above  formula  for  $\pi_n$  in  this  case it follows that
$\pi_n(x_1\otimes\dots\otimes x_{n+2})=0$,  if $x_1\le x_2$ or $\dots$
or     $x_{n+1}\le     x_{n+2}$.     If    $x_1>\dots>x_{n+2}$    then
$\pi_n(x_1\otimes\dots \otimes x_{n+2})\ne 0$ in general.

So $A_\infty$-algebra structure on $A_*=PX$ may be not degenerated. As
it  was  shown  in  [15],  on  the  polynomial  algebra  $P$  with $n$
generators,   $n\ge   3$,   there   really   exists    none    trivial
$A_\infty$-algebra  structures.  In  the case $n=3$ they are in one to
one correspondence with the Hochschild homology $H^3(P;P)\cong P$.

So none  trivial  elements  $y\in  P$   degenerate   none   degenerated
$A_\infty$-algebra structures.  The  corresponding  $A_\infty$-algebra
denote by $\widetilde P$.

In the capacity of $A$ we may take the cobar construction over the bar
construction over $\widetilde P$,  i.e.  $A=FB\widetilde P$.  Then $A$
will be  a  differential  algebra  which  homology  is  isomorpic   to
$\widetilde P$. So $A$ is the desired algebra.

Consider the   question  about  what  additional  conditions  must  be
satisfied for a differential algebra $A$ to be the  $A_\infty$-algebra
structure on $A_*$ degenerated.  To give answer to this question we
introduce the notion of a homotopy trivial Lie algebra.

A Lie algebra $L$ with a multiplication $\mu\colon  L\otimes  L\to  L$
will   be   called  homotopy  trilial  if  there  are  given  mappings
$$\mu_n\colon L^{\otimes n+1}\to L,$$ increasing dimensions by $n$ and
satisfying             the             following             relations
$$d(\mu_{n+1})(x_0\otimes\dots\otimes            x_{n+1})=\sum_{p+q=n}
\mu(\mu_p(x_{i_0}\otimes\dots\otimes      x_{i_p})\otimes\mu_q(x_{j_0}
\otimes\dots\otimes x_{j_q})),$$ where $\mu_0(x)=x$  and  the  sum  is
takene over all shuffles $I=\{i_0,\dots,i_p\}$,  $J=\{j_0,\dots,j_q\}$
of the set $0,\dots,n+1$, such that $I<J$.

Consider examples of homotopy trivial  Lie  algebras.  Let  $A$  be  a
differential algebra. It may be turned to a Lie algebra by introducing
a new multiplication $\mu\colon  A\otimes  A\to  A$,  by  the  formula
$$\mu(x\otimes  y)=x\cdot  y-  y\cdot  x.$$  Suppose  the  algebra $A$
possesses  a  $\cup_1$-product   $\cup_1\colon   A\otimes   A\to   A$,
satisfying        the        relations       $$\gather       d(x\cup_1
y)=d(x)\cup_1y+x\cup_1d(y)+x\cdot   y-y\cdot   x;\\   x\cup_1(y_1\cdot
y_2)=(x\cup_1y_1)\cdot y_2+y_1\cdot(x\cup_1y_2).\endgather$$ Show that
in this case the corresponding Lie algebra will be homotopy trivial.

Define mappings $\mu_n\colon A^{\otimes n+1}\to A$ putting
$$\mu_n(x_1\otimes x_2\otimes\dots\otimes x_{n+1})=x_1\cup_1(x_2\cup_1
(\dots\cup_1x_{n+1})\dots).$$
Direct calculations show that the required relations are satisfied.

Note that   instead   of   the   distributivity   relation   for   the
$\cup_1$-product we may demand the homotopy  distributivity  relation.
Namely, let  besides the operation $\cup_1$ there are given operations
$$\phi_n\colon    A^{\otimes    n+1}\to    A,\quad    n\ge     1,\quad
\phi_1=\cup_1,$$      satisfying      the      relations      $$\align
d(\phi_{n+1})(x\otimes          y_1\otimes\dots\otimes          y_n)&=
y_1\cdot\phi_n(x\otimes         y_2\otimes\dots\otimes         y_n)+\\
&+\sum_{i=1}^{n-1}\phi(x\otimes    y_1\otimes\dots\otimes     y_i\cdot
y_{i+1}       \otimes\dots\otimes       y_n)+\\      &+\phi_n(x\otimes
y_1\otimes\dots\otimes y_{n-1})\cdot y_n.\endalign $$

Then on $A$ we also may define the structure of a homotopy trivial Lie
algebra     putting     $$\align\mu_n(x_0\otimes\dots\otimes    x_n)&=
x_0\cup_1(x_1\cup_1(\dots(x_{n-1}\cup_1                  x_n)\dots)+\\
&+x_0\cup_1(\dots\cup_1\phi_2(x_{n-2}\otimes(x_{n-1}\otimes    x_n-x_n
\otimes  x_{n-1}))\dots)+\\  &+\dots   \\   &+\sum   \phi_n(x_0\otimes
x_{i_1}\otimes\dots\otimes  x_{i_n}),  \endalign $$ where the last sum
is taken over all permutations of the collection $1,\dots,n$.

Direct calculations show that the required relations are satisfied.

From here  it  follows  that  the  cochain  complex  $C^*(X)$   of   a
topological space  $X$  gives us the example of a homotopy trivial Lie
algebra.  Indeed $C^*(X)$ is an $E_\infty$-algebra [4] and  hence  the
Lie algebra structure on $C^*(X)$ will be homotopy trivial.

Another example  of  a homotopy trivial Lie algebra gives us the cobar
construction $FK$ over a Hopf algebra $K$.  Indeed in this case on the
cobar    construction    $FK$    there    is    defined   distributive
$\cup_1$-product.  On  generators  it  is   given   by   the   formula
$$[x]\cup_1[y]=[x\cdot  y].$$  Therefore  the Lie algebra structure on
$FK$ will be homotopy trivial.

{\bf Theorem 15.} {\sl Let $A$ be a differential algebra for which the
corresponding  Lie  algebra  structure  is  homotopy  trivial  and the
homology   $A_*=H_*(A)$    is    the    polynomial    algebra.    Then
$A_\infty$-structure on $A_*$ is degenerated.}

{\bf Proof.} Let $A_*=PX$. Put in order generators $x\in X$ and denote
by $\xi(x)\in A$ their representatives.  Define a  mapping  $\xi\colon
A_*\to             A$            putting            $$\xi(x_{i_1}\dots
x_{i_n})=\xi(x_{i_1})\dots\xi(x_{i_n}),\quad i_1\le\dots\le i_n.$$

Of course in the case when the algebra $A$  is  not  commutative  this
mapping  is not an algebra mapping.  Our task is to add on the mapping
$\xi$ till an $A_\infty$-mapping.

Defing mappings  $\xi_n\colon  A_*^{\otimes  n+1}\to  A$,  putting  on
generators     $$\xi_n(x_{i_0}\otimes\dots\otimes      x_{i_n})=\cases
\mu_n(\xi(x_{i_0})\otimes\dots\otimes\xi(x_{i_n})),&i_0>\dots>i_n,\\
0,&\text{in   other    cases}.\endcases    $$    If  the  elements
$u_0,\dots,u_n\in  A_*$  satisfy  one of the inequalities
$u_0\le u_1,\dots  u_{n-1}\le  u_n$  then  we  put  $$\xi_n(u_0\otimes
u_1\otimes\dots\otimes u_n)=0.$$

Show that on the other elements the mappings $\xi_n$ is determined  by
the     relations     $$\align     d(\xi_{n+1})(u_0\otimes\dots\otimes
u_{n+1})&=\sum_{i=0}^n      \xi_n(u_0\otimes\dots\otimes      u_i\cdot
u_{i+1}\dots\otimes                                        u_{n+1})-\\
&-\sum_{i=0}^n\xi_i(u_0\otimes\dots\otimes                   u_i)\cdot
\xi_{n-i}(u_{i+1}\otimes\dots\otimes u_{n+1}),\endalign$$ meaning that
the family $\{\xi_n\}$ is an $A_\infty$-mapping from $A_*$ to $A$.

The definition  of  $\xi_n(u_0\otimes\dots\otimes  u_n)$  is inductive
over the total numbers of factors of $u_0,\dots,u_n$.

Denote by  $x$  the minimal generator of the elements $u_0,\dots,u_n$.
Suppose that $u_0$ containes $x$.  If $x=u_0$ then from the inequality
$x\le  u_1$  it follows that $$\xi_n(x\otimes\dots\otimes u_n)=0.$$ If
$u_0=x\cdot u'_0$ then  from  the  above  relations  it  follows  that
$$\xi_n(u_0\otimes\dots\otimes
u_n)=\xi(x)\xi_n(u'_0\otimes\dots\otimes  u_n).$$  So  the  value   of
$\xi_n$   is  determined  by  the  value  on  the  element  of  lesser
filtration.

Similary if $u_k$ containes $x$, i.e. $u_k=x\cdot u'_k$ then we have
$$\align
\xi_n(u_0\otimes\dots\otimes x\cdot u'_k\otimes\dots\otimes u_n)&=
\xi_n(u_0\otimes\dots\otimes x\cdot u_{k-1}\otimes u'_k\otimes\dots
\otimes u_n)+\\
&+\xi_k(u_0\otimes\dots\otimes u_{k-1}\otimes x)\cdot\xi_{n-k}(u'_k
\otimes\dots\otimes u_n).\endalign $$
So the value of $\xi_n$ is determined by the  value  on  the  elements
of lesser filtrations and by the value on the element in
which $x$ is contained by $x\cdot u_{k-1}$.

Repeating this procedure we come  to  the  values  of  $\xi_n$  on  the
elements of lesser filtration and the value on the element in which
$x$ is contained by the first factors.

It remains  the  case when $x=u_n$ and elements $u_0,\dots,u_{n-1}$ do
not contain $x$.  In  this  case  we  need  to  consider  the  minimal
generator  contained by $u_0,\dots,u_{n-1}$ and repeat described above
procedure. Thus we obtain that the values of $\xi_n$ are determined by
the given values and values on the elements of lesser filtrations.

For example direct inductive  calculations  show  that  there  is  the
following formula
$$\gather\xi_1(x_{i_1}\dots x_{i_n}\otimes x_{j_1}\dots x_{j_m})=\\
\sum_{k,l}\xi(x_{i_1}\dots x_{i_{k-1}}x_{j_1}\dots x_{j_{l-1}})
\xi_1(x_{i_k}\otimes x_{j_l})\xi(x_{i_{k+1}})\dots\xi(x_{i_n})
\xi(x_{j_{l+1}})\dots\xi(x_{j_m}).\endgather$$

Since the Milnor coalgebra may be obtained as the cohomology of an
$E_\infty$-operad, we will have

{\bf Corollary.} {\sl The Milnor coalgebra $\Cal K$ has  degenerated
$A_\infty$-algebra structure.}

\vskip .5cm
\centerline{7. Functorial homology operations}
\vskip 6pt
Let $\Delta^*=\{\Delta^n\}$  denotes  the  cosimplicial  object of the
category of chain complexes, consisting of the chaing complexes of the
standard $n$-dimensional  simplices.  Let  further $F$ be a functor in
the category  of  chain  complexes   for   which   there   are   given
transformations
$$\Delta^n\otimes F(X)\to  F(\Delta^n\otimes  X),$$   permuting   with
coface  and codegeneracy operators.  Such functor $F$ will be called a
chain functor.

A transformation $\alpha\colon  F'\to  F''$  of  chain  functor  is  a
transformations    of   functors,   commuting   the   diagrams   $$\CD
\Delta^n\otimes F'(X)@>>> F'(\Delta^n\otimes X)\\ @V1\otimes\alpha  VV
@VV\alpha  V\\ \Delta^n\otimes F''(X)@>>> F''(\Delta^n\otimes X)\endCD
$$

Given chain functor $F$ we may consider mappings $$F(f)\colon  F(X)\to
F(Y),$$ induced  not  only  by  chain  mappings  $f\colon  X\to  Y$  of
dimension zero but dimension $n$ also.  Namely given mapping  $f\colon
X\to Y$  of  dimension  $n$  we  represent  as  the restriction of the
mapping $\widetilde   f\colon\Delta^n\otimes   X\to    Y$    on    the
$n$-dimensional generator $u_n\in\Delta^n$. Then the required mapping
$$F(f)\colon F(X)\to  F(Y)$$ of dimension $n$  will be the restriction
of the   composition
$$\Delta^n\otimes F(X)\to      F(\Delta^n\otimes      X)@>F(\widetilde
f)>>F(Y)$$ on the $n$-dimensional generator $u_n\in\Delta^n$.

Given chain functor $F$ denote by $F_*$ the functor,  corresponding to
a  chain   complex   $X$   the   graded   module   of   its   homology
$F_*(X)=H_*(F(X))$.  The  functor  $F_*$  is not only a functor but an
$A_\infty$-functor,  i.e.  there are  functional  homology  operations
which  assigns  to sequences of chain mappings $f^1\colon X^1\to X^2$,
\dots  ,$f^n\colon  X^n\to  X^{n+1}$  the  mapping  $$F_*(f^n,   \dots
,f^1)=H_*(F(f^n),\dots,F(f^1))  \colon  F_*(X^1)\to F_*(X^{n+1}),$$ of
dimension $n-1$.

A transformation $\alpha\colon F'\to F''$ of chain functors induces an
$A_\infty$-transformation of the $A_\infty$-functor      $F'_*$ to the
$A_\infty$-functor $F''_*$.

{\bf Theorem 16.} {\sl Let $F$  be  a  chain  functor.  Then  for  any
sequence of chain mappings
$f^1\colon X^1\to  X^2$,  \dots  ,  $f^n\colon  X^n\to  X^{n+1}$   the
following formula is taken place
$$H_*(F(f^n), \dots  ,  F(f^1))=\sum  (-1)^{\epsilon}F_*(H_*(f^n,\dots
,f^{n_m+1}),\dots ,H_*(f^{n_1},\dots  ,f^1)),$$ where the sum is taken
over $m$ and $n_1,\dots , n_m$ such that $1\le n_1 < \dots < n_m < n$.}

{\bf Proof.}  Let  $X$  be  chain  complex.  We   take   the   mapping
$F_*(X_*)\to  F(X)$  of  choosing  representatives  as the composition
$$\xi(F)\colon F_*(X_*)\to F(X_*),\quad F(\xi)\colon F(X_*)\to F(X).$$
Similary  we take the projection $F(X)\to F_*(X_*)$ as the composition
$$F(\eta)\colon   F(X)\to   F(X_*),\quad    \eta(F)\colon    F(X_*)\to
F_*(X_*).$$  We  take  the  homotopy $H\colon F(X)\to F(X)$ as the sum
$$F(\xi)\circ h(F)\circ F(\eta)+F(h).$$ Substituting these mappings to
the  formula  of  functional homology operation we obtain the required
formula.

Similary there is the following theorem.

{\bf Theorem  17.}  {\sl   Let   $\alpha\colon   F'\to   F''$   be   a
transformation  of  chain  functors.  Then  for  any sequence of chain
mappings $f^1\colon X^1\to X^2$,\dots,  $f^n\colon X^n\to X^{n+1}$ the
following           formula           is          taken          place
$$H_*(F''(f^n),\dots,\alpha(X^i),\dots,F'(f^1))=H_*(F''(f^n_*),\dots,
\alpha(X^i_*),\dots,F'(f^1_*)).$$}

A functor $F$ will be called formal if homotopies $h(F)$ may be chosen
in such a way that for any mapping of graded  modules  $f\colon  M'\to
M''$  the following relation is satisfied $$h''(F)\circ F(f)=F(f)\circ
h'(F).$$

Note that from the last relation  it  follows  $$\eta''(F)\circ  F(f)=
F_*(f)\circ \eta'(F),~ F(f)\circ\xi'(F)=\xi''(F)\circ F_*(f).$$

From the definition of functorial homology operations directly follows
that  if  $F$ is a formal functor restricted on the category of graded
modules then the $A_\infty$-structure on $F_*$ is degenerated. In this
case             there            is            the            formula
$$F_*(f^n,\dots,f^1)=F_*(H_*(f^n,\dots,f^1)).$$

A transformation $\alpha\colon F'\to F''$ of chain  functors  will  be
called formal if homotopies $h(F')$ and $h(F'')$ may be chosen in such
a way that the following relation  is  satisfied  $$h''(F)\circ\alpha=
\alpha\circ h'(F).$$

If $\alpha\colon  F'\to  F''$  is  a  formal  transformation  then  the
$A_\infty$-transformation structure   from   $F'_*$   to   $F''_*$  is
degenerated.

In this   case   there   is   the   formula   $$\alpha_*(X_{n+1})\circ
F'_*(f^n,\dots,f^1)= F''_*(f^n,\dots,f^1)\circ\alpha_*(X_1).$$

\vskip .5cm
\centerline{8. Homology operations for the operad $E$}
\vskip 6pt
Show that the functors $\underline E$,  $\overline E$ corresponding to
an  $E_\infty$-operad $E$ are chain.  To do it we define the family of
mappings     $$\Delta^n\otimes\underline     E(j;X)\to      \underline
E(j;\Delta^n\otimes    X)$$   to   be   the   compositions   $$\gather
\Delta^n\otimes\underline                       E(j;X)=\Delta^n\otimes
E(j)\otimes_{\Sigma_j}    X^{\otimes   j}@>1\otimes\nabla\otimes   1>>
\Delta^n\otimes  E(j)\otimes   E(j)\otimes_{\Sigma_j}X^{\otimes   j}\\
@>\tau\otimes   1\otimes  1>>\Delta^{n\otimes  j}\otimes  E(j)\otimes_
{\Sigma_j}X^{\otimes  j}\to  E(j)\otimes_{\Sigma_j}(\Delta^n   \otimes
X)^{\otimes j}=\underline E(j;\Delta^n\otimes X),  \endgather $$ where
$\tau\colon\Delta^n\otimes   E(j)\to\Delta^{n\otimes   j}$    is    an
$E$-coalgebra structure on the complex $\Delta^n$.

     Direct verification  show  that  the   required   relations   are
satisfied.

     Similary define mappings
$$\Delta^n\otimes\overline  E(j;X)\to \overline E(j;\Delta^n\otimes X)$$
or, that is the same, mappings
$$E(j)\otimes\Delta^n\otimes Hom_{\Sigma_j}(E(j);X^{\otimes     j})\to
(\Delta^n\otimes    X)^{\otimes    j}$$ to be the compositions
$$\gather
E(j)\otimes\Delta^n\otimes     Hom_{\Sigma_j}(E(j);X^{\otimes      j})
@>\nabla\otimes   1\otimes  1>>E(j)\otimes  E(j)\otimes\Delta^n\otimes
Hom_{\Sigma_j}(E(j);X^{\otimes j})\\  \to  E(j)\otimes\Delta^{n\otimes
j}\otimes   Hom_{\Sigma_j}(E(j);X^   {\otimes   j})\to(\Delta^n\otimes
X)^{\otimes j}. \endgather $$

     Direct verification  show  that  the   required   relations   are
satisfied. So we have

{\bf Theorem  18.}  {\sl The functors $\underline E_*$, $\overline E_*$
are $A_\infty$-functors.}

     Our aim is to calculate the functional  homology  operations  for
the functors $\underline E_*$,  $\overline E_*$. It means that for any
sequence  $$f^1:X^1\to  X^2,\dots  ,f^n:X^n\to  X^{n+1}$$   of   chain
mappings  we need to calculate the mappings $$\underline E_*(f^n,\dots
,f^1)\colon\underline              E_*(X^1)\to              \underline
E_*(X^{n+1}),\quad\overline     E_*(f^n,\dots     ,f^1)\colon\overline
E_*(X^1)\to \overline E_*(X^{n+1}).$$

Consider firstly the functor  $\underline  E(2;  -  )$ correponding
to a complex $X$  the complex   $$\underline      E(2;X)      =
E(2)\otimes_{\Sigma_2}X\otimes  X,$$ where  $E(2)$ -- $\Sigma_2$-free
and acyclic  complex  with  generators  $e_i$  of dimensions $i$. A
differential is defined by the formula $$d(e_i)=e_{i-1}+e_{i-1}T,\quad
T\in\Sigma_2.$$

The homology $\underline E_*(2; - )$ of this functor, as it was poined
out above, is not only a functor but an $A_{\infty}$-functor. It means
that  for  any  sequence  of  chain  mappings  $$f^1:X^1\to  X^2,\dots
,f^n:X^n\to   X^{n+1}$$   there   is   the   operation    $$\underline
E_*(2;f^n,\dots    ,f^1)\colon\underline    E_*(2;X^1)\to   \underline
E_*(2;X^{n+1}).$$ Let us calculate these operations.

Note that for a chain complex $X$ there is an isomorphism $$\underline
E_*(2,X)\cong\underline E_*(2,X_*).$$ If $X_*$ is a graded module then
$\underline E_*(2;X_*)$ is the sum of two factors. The first factor is
the  quotient  module $X_*\cdot X_*$ of the tensor product $X_*\otimes
X_*$ up to permutation of factors.  The second factor  is  the  module
generated  by  the  elements of the form $e_i\times y_n$,  $i\ge 1$ of
dimensions $i+2n$. The elements $y_n\cdot y_n\in X_*\cdot X_*$ will be
also denoted as $e_0\times y_n$.

Let $\xi\colon X_*\to X$,  $\eta\colon X\to X_*$, $h\colon X\to X$ are
mappings giving a chain equivalence between  $X$ and  $X_*$.
Denote by  $$E(\xi)\colon  E(2,X_*)\to  E(2,X),~
E(\eta)\colon E(2,X)\to E(2,X_*),~ E(h)\colon   E(2,X)\to   E(2,X)$$
the mappings defined by the formulas
$$\gather E(\xi)(e_i\otimes y_1\otimes  y_2)=e_i\otimes\xi(y_1)\otimes
\xi(y_2),~ E(\eta)(e_i\otimes x_1\otimes x_2)=e_i\otimes\eta(x_1)
\otimes \eta(x_2),\\ E(h)(e_i\otimes x_1\otimes  x_2)=e_i\otimes (x_1
\otimes h(x_2)+h(x_1)\otimes    \xi\eta(x_2))+e_{i-1}\otimes    h(x_1)
\otimes h(x_2).\endgather $$
It is clear they give a chain equivalence $\underline
E(2,X)\simeq\underline E(2,X_*)$.

Define mappings     $$\gather\xi(E)\colon\underline      E_*(2;X_*)\to
\underline  E(2;X_*),~  \eta(E)\colon\underline  E(2;X_*)\to\underline
E_*(2;X_*),\\        h(E)\colon\underline        E(2;X_*)\to\underline
E(2;X_*).\endgather$$  To do it firstly we we choose an ordering basis
$\{y\}$  in  $X_*$.  Then  define   the   mapping   $\xi(E)$   putting
$$\xi(E)(e_i\times  y)=e_i\otimes  y\otimes y;~ \xi(E) (y_1\cdot y_2)=
e_0\otimes  (y_1\otimes  y_2),~y_1\le  y_2.$$   Define   the   mapping
$\eta(E)$   putting  $$\eta(E)(e_i\otimes  y_1\otimes  y_2)  =  \cases
e_i\times y_1,  & y_1=y_2 \\  y_1\cdot  y_2,  &  y_1<y_2,  i=0  \\  0,
&\text{in other cases} \endcases $$

Define the mapping     $h(E)$ putting   $$h(E)(e_i\otimes
y_1\otimes y_2)=\cases e_{i+1}\otimes y_2\otimes y_1,  & y_1>y_2,\\ 0,
& \text{in other cases}\endcases $$

     Direct calculations  show  that  the   required   relations   are
satisfied.

The mappings    $$\gather    E(\xi)\circ\xi(E)\colon     E_*(2,X_*)\to
E(2,X),\quad   \eta(E)\circ   E(\eta)\colon   E(2,X)\to  E_*(2,X_*),\\
E(\xi)\circ h(E)\circ E(\eta)+E(h)\colon E(2,X)\to E(2,X)\endgather $$
give us a chain equivalence between $E(2,X)$ and $E_*(2,X_*)$.

From the general formula of functional homology operations for a chain
functor it  follows  that for the functor $\underline E(2,-)$ the
following formula is taken place
$$\underline E_*(2,f^n,\dots,f^1)=\sum                      \underline
E_*(2,H_*(f^n,\dots,f^{n_m+1}),\dots ,H_*(f^{n_1},\dots      ,f^1)),$$
where the  sum  is taken over all $m$ and $n_1,\dots ,  n_m$ such that
$1\le n_1 < \dots < n_m < n$.

For a graded module $X$ with fixed  ordering  basis  $\{x_i\}$  define
mappings $p\colon  X\otimes X\to X$, $q\colon
X\to X\otimes  X$,  $r\colon  X\otimes  X\to  X\otimes   X$ putting
$q(x_i)=x_i\otimes x_i$ and
$$p(x_i\otimes x_j)=\cases x_i,&i=j,\\0,&i\ne j,\endcases\quad
r(x_i\otimes x_j)=\cases x_j\otimes x_i,&i>j,\\0,&i\ge j.\endcases $$

For a sequence of mappings   $f^1\colon   X^1\to  X^2,\dots,
f^n\colon X^n\to X^{n+1}$ of  graded  modules  with  ordering  basises
define the  mapping  $(f^n,\dots,f^1)\colon  X^1\to X^{n+1}$ putting
$$(f^n,\dots,f^1)=p\circ(f^n)^{\otimes            2}             \circ
r\circ(f^{n-1})^{\otimes   2}\circ\dots\circ   r\circ   (f^1)^{\otimes
2}\circ q.$$

Directly from the definition of the homology operations it follows

{\bf Theorem  19.}  {\sl  For a sequence of mappings $f^1\colon X^1\to
X^2,\dots,f^n\colon X^n\to X^{n+1}$ of graded  modules  the  following
formula  is  taken  place $$\underline E_*(2;f^n,\dots ,f^1)(e_i\times
x)= e_{i+n-1}\times (f^n,\dots ,f^1)(x).$$}

To obtain the correpondig formula for the functor $\underline E_*$  it
needs to use a monad structure
$\underline\gamma_*\colon\underline E_*\circ\underline          E_*\to
\underline E_*$ and the formula
$$\underline E_*(f^n,\dots,f^1)\circ\gamma_*=\sum           \underline
\gamma_*\circ\underline E_*(\underline  E_*(f^n,\dots,f^{n_m+1}),\dots
,\underline E_*(f^{n_1},\dots ,f^1)),$$ where the sum  is  taken  over
all $m$ and $n_1,\dots , n_m$ such that $1\le n_1 < \dots < n_m < n$.

Passing to  the  Dyer-Lashof algebra $\Cal R$ we obtain the operations
$$\Cal   R(f^n,\dots,f^1)\colon\Cal   R\times    X^1\to\Cal    R\times
X^{n+1},$$ which on the generators $Q^i$ are expressed by the formulas
$$\Cal  R(f^n,\dots,f^1)(Q^i\otimes  x)=  Q^{i+n-1}\otimes  (f^n,\dots
,f^1)(x).$$

Dually for the functor $\overline   E_*(2;-)$ there is

{\bf Theorem   20.}   {\sl For a sequence of mappings
$f^1\colon X^1\to  X^2,\dots,f^n\colon  X^n\to  X^{n+1}$   of   graded
modules the following formula is taken place
$$\overline  E_*(2;f^n,\dots  ,f^1)(\overline  e_i\times  x)=
\cases\overline  e_{i-n+1}\times  (f^n,\dots  ,f^1)(x),&  i\ge   n-1\\
0,&\text{in other cases} \endcases $$}

To obtain the corresponding formula for the functor $\overline E_*$ it
needs to use a comonad structure
$\overline\gamma_*\colon\overline E_*\to\overline    E_*\circ\overline
E_*$  and the formula
$$\overline\gamma_*\circ\overline E_*(f^n,\dots,f^1)=
\sum \overline  E_*(\overline E_*(f^n,\dots,f^{n_m+1}),\dots,\overline
E_*(f^{n_1},\dots ,f^1))\circ\overline\gamma_*,$$  where  the  sun  is
taken over all $m$ and $n_1,\dots , n_m$ such that $1\le n_1 < \dots <
n_m < n$.

Passing to the Milnor coalgebra  $\Cal  K$ we obtain the operations
$$\Cal K(f^n,\dots,f^1)\colon\Cal     K\times    X^1\to\Cal    K\times
X^{n+1},$$ which are expressed by the formulas
$$\Cal K(f^n,\dots,f^1)(y\otimes      x)=     y\cdot\xi_1^{n-1}\otimes
(f^n,\dots ,f^1)(x).$$

Consider operations associated with a comultiplication  $\nabla$  of
the   Milnor  coalgebra  $\Cal  K$.  Denote  $$\nabla(n)=\nabla\otimes
1\dots\otimes         1-\dots+(-1)^{n-1}1\otimes\dots          \otimes
1\otimes\nabla\colon \Cal K^{\times n}\to\Cal K^{\times n+1}.$$

Direct calculations        show        that       the       operations
$$(\nabla(n),\dots,\nabla)\colon \Cal  K\to\Cal  K^{\times  n+1},~n\ge
2,$$  are trivial on the elements $\xi_i^{2^k}$.  However on the other
elements these operations in general  are  not  trivial.  For  example
there     is     the     formula     $$(\nabla(2),\nabla)(\xi_i\xi_j)=
\xi_{j-i}^{2^i}\xi_0^{2^i}\otimes\xi_i\xi_0^{2^i}\otimes\xi_i,~i<j.$$

Denote by $\widetilde\nabla$ a comultiplication in the tensor  product
$\Cal    K\otimes\Cal    K$,   $$\widetilde\nabla=(1\otimes   T\otimes
1)(\nabla\otimes\nabla).$$                                         Put
$$\widetilde\nabla(n)=\widetilde\nabla\otimes  1\dots\otimes  1-\dots+
(-1)^{n-1}1\otimes\dots\otimes  1\otimes\widetilde\nabla\colon   (\Cal
K\otimes\Cal  K)^{\times  n}\to(\Cal  K\otimes\Cal  K)^{\times n+1}.$$
Consider           the            operations            $$(\pi^{\times
n+1},\widetilde\nabla(n),\dots,\widetilde\nabla)\colon \Cal K^{\otimes
2}\to \Cal K^{\times n+1}.$$ Its restriction on the elements $x\otimes
x\in  \Cal  K\otimes  \Cal K$ we denote by $$\Psi^n\colon\Cal K\to\Cal
K^{\times n+1}.$$

From the  formula  of  a  comultiplication  in  the  Milnor  coalgebra
directly        follows       the       formula       $$\Psi^1(\xi_n)=
\sum_{i<j}\xi_{n-i}^{2^i}\xi_{n-j}^{2^j}\otimes\xi_i\xi_j,$$   or   in
more        general        case       $$\Psi^1(\xi_n^{2^m})=\sum_{i<j}
\xi_{n-i}^{2^{i+m}}\xi_{n-j}^{2^{j+m}}\otimes\xi_i^{2^m}\xi_j^{2^m}.$$
In  particular  for  the  primitive elements $\xi_1^{2^m}\in\Cal K$ we
have                            the                            formula
$$\Psi^1(\xi_1^{2^m})=\xi_1^{2^m}\xi_0^{2^m}\otimes\xi_1^{2^m}.$$

Similary for the operation $\Psi^2$ we have the formula
$$\Psi^2(\xi_n^{2^m})=\sum\Sb i<j\\ k>l\endSb
\xi_{n-i}^{2^{i+m}}\xi_{n-j}^{2^{j+m}}
\otimes\xi_{i-k}^{2^{k+m}}\xi_{j-l}^{2^{l+m}}
\otimes\xi_k^{2^m}\xi_l^{2^m}.$$
In particular for the primitive  elements  $\xi_1^{2^m}\in\Cal  K$  we
have the formula $$\Psi^2(\xi_1^{2^m})=0.$$ And so on.

\vskip .5cm
\centerline{9. $\cup_\infty-A_\infty$-Hopf algebras}
\vskip 6pt
To calculate  higher  differentials  of the Adams spectral sequence we
need to use not only the action of the Dyer-Lashof  algebra,  but  the
$E_\infty$-structure. However  this structure is too complicate.  Some
of the calculations were made in [6].  Here we'll use only a  part  of
the   $E_\infty$-structure  consisting  of  $\cup_i$-products.

A chain  complex  $A$  will be called a $\cup_\infty$-algebra if there
are given operations $\cup_i\colon A\otimes A\to A$,  $i\ge 0$, called
$\cup_i$-products,  increasing  dimensions  by  $i$ and satisfying the
relation
$$d(x\cup_iy)=d(x)\cup_iy+x\cup_id(y)+x\cup_{i-1}y+y\cup_{i-1}x.$$

A differential  coalgebra  $K$  will  be  called  a $\cup_\infty$-Hopf
algebra if there are given  $\cup_i$-products  $\cup_i\colon  K\otimes
K\to      K$      satisfying      the      distributivity     relation
$$\nabla(x\cup_iy)=\sum_k(x'\cup_{i-k}T^ky')\otimes(x''\cup_ky''),$$
where $\nabla(x)=\sum x'\otimes x''$, $\nabla(y)=\sum y'\otimes y''$,
$T\colon K\otimes K\to K\otimes K$ is the permutation mapping, $T^k$
it's $k$-th iteration.

{\bf Theorem   21.}   {\sl   The   cobar   construction  $FK$  over  a
$\cup_\infty$-Hopf algebra $K$ is  a  $\cup_\infty$-algebra.  Moreover
$\cup_i$-products   $\cup_i\colon   FK\otimes   FK\to   FK$   uniquily
determined by the formula
$$[x]\cup_i[y]=\cases  [x\cup_{i-1}y],&i\ge  1,\\  [x,y],&i=0.
\endcases $$
and the relations
$$\align
(x_1x_2)\cup_i[y]&=(x_1\cup_i[y])x_2+x_1(x_2\cup_i[y]),\\
(x_1x_2)\cup_i(y_1y_2)&=\sum_k(x_1\cup_{i-k}T^ky_1)
(x_2\cup_ky_2)+\\
&+(x_1\cup_i(y_1y_2))x_2+x_1(x_2\cup_i(y_1y_2))+\\
&+((x_1x_2)\cup_iy_1)y_2+y_1((x_1x_2)\cup_iy_2)+\\
&+x_1(x_2\cup_iy_1)y_2+y_1(x_1\cup_iy_2)x_2,\endalign $$
where $x_1,x_2,y_1,y_2\in FK$, $y\in K$, $i\ge 1$.}

Indeed, the products $[x_1,\dots,x_n]\cup_i[y]$ are determined by  the
first relation  $$[x_1,\dots,x_n]\cup_i[y]=\sum_{k=1}^n[x_1,\dots,x_k
\cup_iy,\dots,x_n]$$

From the second relation it follows that to define $\cup_i$-products
in general case,  i.e.  $$[x_1,\dots,x_n]\cup_i[y_1,\dots,y_m]$$ it is
sufficient        to        define        only       $\cup_i$-products
$[x]\cup_i[y_1,\dots,y_m]$.

We have
$$\align d([y_1,\dots,y_m]\cup_{i+1}[x])&=
\sum_{k=1}^m[y_1,\dots,d(y_k),\dots,y_m]\cup_{i+1}[x]+\\&+
\sum_{k=1}^m[y_1,\dots,y'_k,y''_k,\dots,y_m]\cup_{i+1}[x]+\\&+
[y_1,\dots,y_m]\cup_{i+1}([d(x)]+[x',x''])+\\&+
[y_1,\dots,y_m]\cup_i[x]+[x]\cup_i[y_1,\dots,y_m].\endalign $$

Thus the   product  $[x]\cup_i[y_1,\dots,y_m]$  is  expressed  through
already defined   products   and   products  of  the  elements  lesser
dimensions. Hence $\cup_i$-products are determined by induction.

So this  theorem  gives  us  the  formulas for $\cup_i$-product in the
cobar construction.  However they are inductive and not so simple even
in the case when  higher  $\cup_i$-products  ($i\ge  1$)  on  $K$  are
trivial, i.e. when $K$ is a commutative Hopf algebra.

An $\cup_\infty$-Hopf algebra $K$ will be called commutative if the
coproduct $\nabla\colon K\to K\otimes K$ is comutative.

{\bf Theorem 22.} {\sl The cobar construction $FK$ over a commutative
$\cup_\infty$-Hopf algebra is a commutative $\cup_\infty$-Hopf algebra.
So the cobar construction over a commutative $\cup_\infty$-Hopf algebra
may be iterated.}

{\bf Proof.} Define the coproduct $\nabla\colon FK\to FK\otimes FK$ putting
$$\nabla[x_1,\dots,x_n]=\sum[x_{i_1},\dots,x_{i_p}]\otimes [x_{j_1},
\dots,x_{j_q}],$$ where the sum is taken over all $(p,q)$-shuffles of
$1,2,\dots,n$. Direct calculatons show that the required relations are
satisfied.

Consider the   question   about  the  structure  on  the  homology  of
a $\cup_\infty$-Hopf algebra.

Consider the  question  about  the  structure  on  the  homology  of a
$\cup_\infty$-Hopf algebra.  It is clear that on the homology $K_*$ of
a  $\cup_\infty$-Hopf  algebra  $K$  there  are  $\cup_\infty$-algebra
structure,  consisting  of  the  operations  $$\cup_i\colon  K_*\otimes
K_*\to  K_*$$  and  $A_\infty$-coalgebra  structure,  consisting of the
operations $$\nabla_n\colon K_*\to K_*^{\otimes n+2}.$$ But besides that
there  are another operations of the form $$\Psi_{i,n}\colon K_*\otimes
K_*\to K_*^{\otimes n+2}.$$

To   describe these operations and relations between them we introduce
the notion of a $\cup_\infty-A_\infty$-Hopf algebra.

An $A_\infty$-coalgebra $K$      will      be      called      an
$\cup_\infty-A_\infty$-Hopf algebra     if on the cobar construction
$\widetilde FK$ there is given $\cup_\infty$-algebra structure
satisfying the relations
$$\align
(x_1x_2)\cup_i[y]&=(x_1\cup_i[y])x_2+x_1(x_2\cup_i[y]),\\
(x_1x_2)\cup_i(y_1y_2)&=\sum_k(x_1\cup_{i-k}T^ky_1)
(x_2\cup_ky_2)+\\
&+(x_1\cup_i(y_1y_2))x_2+x_1(x_2\cup_i(y_1y_2))+\\
&+((x_1x_2)\cup_iy_1)y_2+y_1((x_1x_2)\cup_iy_2)+\\
&+x_1(x_2\cup_iy_1)y_2+y_1(x_1\cup_iy_2)x_2,\endalign $$
where $x_1,x_2,y_1,y_2\in FK$, $y\in K$, $i\ge 1$.

{\bf Theorem 23.} {\sl If $K$ is a $\cup_\infty$-Hopf algebra then its
homology $K_*=H_*(K)$ is $\cup_\infty-A_\infty$-Hopf algebra and there
is  an  equivalence  of  $\cup_\infty$-algebras $\widetilde FK_*\simeq
FK$.}

{\bf Proof.} It is known [14] that the homology $K_*$ of a differential
coalgebra $K$ is $A_\infty$-coalgebra and there are algebra mappings
$\xi\colon \widetilde FK_*\to FK$, $\eta\colon FK\to\widetilde FK_*$ and
an algebra chain homotopy $h\colon FK\to FK$ such that
$\eta\circ\xi=Id$, $d(h)=\xi\circ\eta-Id$. So we need to define on
$\widetilde FK_*$ the $\cup_i$-products. Put on generators
$$[x]\cup_i[y]=\eta(\xi[x]\cup_i\xi[y]).$$
On the other elements the $\cup_i$-products determines by the relations.

Applying this theorem to the Milnor coalgebra we obtain

{\bf Theorem  24.}  {\sl  The  Milnor  coalgebra  $\Cal  K$   possesses
$\cup_\infty-A_\infty$-Hopf algebra structure.  The  homology  of  the
corresponding cobar construction $\widetilde F\Cal K$ is isomorphic to
the  $E^\infty$  term  of  the   Adams   spectral   sequence.}

\vskip .5cm
\centerline{10. The calculation of the differentials}
\vskip 6pt
Apply developed methods to calculate higher differentials of the Adams
spectral  sequence or,  that is the same to calculate the differential
in $\widetilde F\Cal K$.

Recall that there are the following formulas for the $\cup_i$-products
on  the  Milnor coalgebra
$$   \xi_n\cup_i\xi_n   =   \cases   \xi_{n+k}\xi_n,   &   i
=2^{n+k}-2^n;\\ \xi_{n+k}\xi_{n-1}, & i = 2^{n+k}-2^n-2^{n-1};\\ \dots
& \dots \\ \xi_{n+k}\xi_0,& i = 2^{n+k}-2^n-\dots -1;\\  0,&  \text{in
other cases.} \endcases $$

Thus any element of  $\Cal  K$  may  be  obtained  from  $\xi_0$  by
applying $\cup_i$-products. Hence we have

{\bf Theorem 25.} {\sl  The  formulas  for  $\cup_i$-products  in  the
Milnor coalgebra and the relations for $\cup_\infty$-algebra structure
in the cobar construction $\widetilde F\Cal  K$  completely  determine
the differential in $\widetilde F\Cal K$.}

However the  formulas  for  the  differential  are  inductive and very
complicated.  So the next step in the calculation of the  differential
is to replace the Milnor coalgebra $\Cal K$ and the cobar construction
$\widetilde F\Cal K$ by more simply objects.

To do it consider the filtration of $\widetilde F\Cal K$  putting  the
filtration  of  the elements $\xi_{i_1}\cdot\dots\cdot \xi_{i_n}\in K$
to be equal $n$.  Then the first term of  the  corresponding  spectral
sequence will  be isomorphic to the polynomial algebra $PS^{-1}X$ over
the module $X$ generated by the elements $\xi_i^{2^k}$.  We'll  denote
elements of  $PS^{-1}X$  as  elements of the cobar construction,  i.e.
$[x_1,\dots,x_n]$, $x_i\in X$.

Note that there is an algebra mapping $\eta\colon F\Cal K\to PS^{-1}X$,
given       by       the       formula
$$\eta [x]=\cases   [\xi_i^{2^k}],&x=\xi_i^{2^k},\\0,&\text{in   other
cases.}\endcases $$

The inverse  mapping  $\xi\colon  PS^{-1}X\to \Cal FK$ may be given by
the formula
$$\xi([x_1,\dots,x_n])=[x_1,\dots,x_n],\quad x_1\le\dots\le x_n.$$
It is not an  algebra  mapping,  but  if  $x\le  y$  then  $\xi(x\cdot
y)=\xi(x)\cdot\xi(y)$.

From here and from the  Perturbation  Theory it follows

{\bf Theorem 26.} {\sl The  polynomial  algebra  $PS^{-1}X$  possesses
a $\cup_\infty$-algebra structure determined on the generators by  the
formulas
$$[x]\cup_i[y]=\cases 0,&i\ge 1,~x<y,\\
[x\cup_{i-1}x],&i\ge 1,~x=y,\\
[x,y],&i=0,~x\le y\endcases $$
and satisfying the relations
$$\align
(x_1x_2)\cup_i[y]&=(x_1\cup_i[y])x_2+x_1(x_2\cup_i[y]),\\
(x_1x_2)\cup_i(y_1y_2)&=\sum_k(x_1\cup_{i-k}T^ky_1)
(x_2\cup_ky_2)+\\
&+(x_1\cup_i(y_1y_2))x_2+x_1(x_2\cup_i(y_1y_2))+\\
&+((x_1x_2)\cup_iy_1)y_2+y_1((x_1x_2)\cup_iy_2)+\\
&+x_1(x_2\cup_iy_1)y_2+y_1(x_1\cup_iy_2)x_2,\endalign $$
where $x_1,x_2,y_1,y_2\in PS^{-1}X$, $y\in X$, $i\ge 1$.

The homology of the corresponding  complex  $\widetilde  PS^{-1}X$  is
isomorphic  to  the  homology of $\widetilde F\Cal K$ and hence to the
$E^\infty$ term of the Adams spectral sequence.}

This theorem gives us the inductive formulas for the differential  and
$\cup_i$-products  in  $\widetilde  PS^{-1}X$.

Since any   element  of  $\widetilde  PS^{-1}X$  may  be  obtained  from
$[\xi_0]$ by  applying  $\cup_i$-products,  we have

{\bf Theorem 27.} {\sl  The  formulas  for  $\cup_i$-products  in  the
module  $X$  and  the  relations  for $\cup_i$-products in $\widetilde
PS^{-1}X$  completely  determine  the  differential   in   $\widetilde
PS^{-1}X$.}

Denote $[\xi_n^{2^m}]$  by  $h_{n,m}$.  Note  that  in  Adams notation
$h_{1,m}=h_m$.

Using the formula $h_{0,0}\cup_{2^n}h_{0,0}=h_{n,0}$ we obtain

{\bf Theorem 28.} {\sl For the differential in  $\widetilde  PS^{-1}X$
on    the    elements    $h_{n,0}$   there   is   the   next   formula
$$d(h_{n,0})=\sum_{i=1}^{n-1}h_{n-i,i}h_{i,0},$$  i.e.  there  are  no
higher differentials on the elements $h_{n,0}$.}

By induction,  using the formula $h_{n,m}\cup_1h_{n,m}=h_{n,m+1}$,  we
obtain the following theorem.

{\bf Theorem 29.} {\sl For the differential in  $\widetilde  PS^{-1}X$
on  the  elements  $h_{n,m}$,  $m\ge  1$,  there  is  the next formula
$$\align
d(h_{n,m})=&\sum_{i=1}^{n-1}h_{n-i,m+i}h_{i,m}+
\sum_{i=1}^{m-1}h_{i,0}h_{n,m-i}^{2^i}+
\sum_{i=1}^{m-2}h_{1,i-1}^2h_{n,m-i}^{2^i}+\\
+&\sum_{i=2}^{m-2}h_{1,0}^{2^i}h_{n,m-1}^2+
\sum_{i=1}^{n-1}h_{n-i,i}^{2^m}h_{i+m,0}+\\
+&\sum_{i=1}^{n-1}h_{n-i,m+i-1}^2h_{i+1,0}^{2^{m-1}}\quad (m>1)+\\
+&\sum_{i=1}^{n-1}h_{n-2,i+1}^{2^{m-1}}h_{i+1,m-2}^2\quad (m>2).
\endalign $$}

In particular,  for  the  Adams  elements  $h_m=h_{1,m}$  we  have the
following formula
$$\align
d(h_m)&=\sum_{i=1}^{m-1}h_{i,0}h_{m-i}^{2^i}+
\sum_{i=1}^{m-2}h_{i-1}^2h_{m-i}^{2^i}+
\sum_{i=2}^{m-2}h_0^{2^i}h_{m-1}^2=\\
&=h_0h_{m-1}^2+\text{elements of greater filtration}
\endalign $$

The first member $h_0h_{m-1}^2$ corresponds to the second differential
of the Adams spectral sequence. So we have the Adams formula
$$d_2(h_m)=h_0h_{m-1}^2.$$

Of course it may  be  obtained  directly  without  using  the  general
formula. Namely,
$$\align d(h_m)&=d(h_{m-1}\cup_1h_{m-1})=\\
&=(h_{-1}h_{m-1}+\dots)\cup_2(h_{-1}h_{m-1}+\dots)=\\
&=(h_{-1}\cup_2h_{-1})h_{m-1}^2+\dots=\\
&=h_0h_{m-1}^2+\dots\endalign $$
where dots denote the elements of greater filtration.
Hence $d_2(h_m)=h_0h_{m-1}^2$.

\vskip .5cm
\centerline{REFERENCES}
\vskip 6pt
1. Adams J.F.  On the  structure  and  applications  of  the  Steenrod
algebra. Comm. Math. Helv. 1958, v.32, p.180--214.

2. Stasheff J.D.  Homotopy associativity of $H$-spaces.  Trans. Amer.
Math. Soc. 1963, v. 108, N 2, p.275--312.

3. May J.P.  The geometry of iterated loop  spaces.  Lect.  Notes  in
Math. 1972, v.271.

4. Smirnov V.A.  On the cochain complex of a  topological  space.  Mat.
Sbornik, 1981, v.115, N.1, p.146--158.

5. Smirnov V.A. Homotopy theory of coalgebras. Izvestia AN SSSR,
1985, v.49, N.6, p.1302--1321.

6. Smirnov V.A.  Secondary operations in the homology of  the  operad
$E$. Izvestia RAN, 1992, v.56, N 2, p.449-468.

7. Steenrod  N.E.  Cohomology  invariants of mappings.  Ann.  of Math.
1949, v.50, p.954--988.

8. Peterson F.P.  Functional cohomology operations. Trans. Amer. Math.
Soc. 1957, v.86, p.187--197.

9. Smirnov V.A. Functional homology operations and weak homotopy type.
Mat. zametki, 1989, v.45, N.5, p.76--86.

10. Bousfield A.K., Kan D.M. The homotopy spectral sequence of a spaces
with coefficients in a ring. Topology, 1972, v.11, p.79--106.

11. Gugenheim V.K.,  Lambe L.A.,  Stasheff J.D. Perturbation theory in
Differen\-tial Homological Algebra. Ill. J. of Math. 1991, v. 35, N 3,
p.357--373.

12. Araki S., Kudo T. Topology  of  $H_n$-spaces  and  $H_n$-squaring
operations. Mem.  Fac.  Sci.  Kyusyu Univ.,  Ser.  A, 1956, v.10, N 2,
p.85--120.

13. Dyer E., Lashof R.K. Homology of iterated loop spaces. Amer. Jour.
of Math. 1962, v.84, N 1, p.35--88.

14. Kadeishvili T.V. On the homology theory of fibre spaces.
UMN. 1980. v. 35, N. 3. p. 183--188.

15. Smirnov  V.A.  $A_\infty$-structures  and  the  functor  $\Cal  D$.
Izvestiya RAN, 2000, N 5 ... .

\vskip 1cm
E-mail: V.Smirnov\@g23.relcom.ru

\enddocument